\documentclass[11pt]{amsart}
\usepackage{amsmath,amsfonts,amsthm,amssymb,
amsthm,graphicx,mathtools,amscd}
\usepackage{hyperref}

\usepackage[mathscr]{eucal}
\usepackage{graphics}
\usepackage{color}
\usepackage{epsfig}

\usepackage[nocompress]{cite}

\DeclareSymbolFont{iwonaletters}{OML}{iwona}{m}{it}
\DeclareMathSymbol{\bdel}{\mathalpha}{iwonaletters}{"E}

\newtheorem{theorem}{Theorem}[section]
\newtheorem{lemma}[theorem]{Lemma}

\newtheorem{proposition}[theorem]{Proposition}

\newtheorem{remark}[theorem]{Remark}

\newcommand{\beginsec}{
\setcounter{equation}{0}
}

\newcommand{\la}{\lambda}
\newcommand{\eps}{\varepsilon}

\newcommand{\al}{\alpha}

\newcommand{\gam}{\gamma}

\newcommand{\sig}{\sigma}
\newcommand{\del}{\delta}
\newcommand{\om}{\omega}
\newcommand{\Gam}{\mathnormal{\Gamma}}

\newcommand{\Om}{\mathnormal{\Omega}}

\newcommand{\N}{{\mathbb N}}

\newcommand{\R}{{\mathbb R}}

\newcommand{\Z}{{\mathbb Z}}

\newcommand{\EE}{{\mathbb E}}
\newcommand{\E}{{\mathbb E}}

\newcommand{\PP}{{\mathbb P}}

\newcommand{\calF}{{\mathcal F}}

\newcommand{\calI}{{\mathcal I}}

\newcommand{\calN}{{\mathcal N}}

\newcommand{\calP}{{\mathcal P}}
\newcommand{\calR}{{\mathcal R}}
\newcommand{\calS}{{\mathcal S}}

\newcommand{\calZ}{{\mathcal Z}}

\newcommand{\skp}{\vspace{\baselineskip}}

\newcommand{\diag}{{\rm diag}}

\newcommand{\To}{\Rightarrow}

\newcommand{\iy}{\infty}
\newcommand{\up}{\uparrow}

\newcommand{\cadlag}{c\`adl\`ag }

\newcommand{\noi}{\noindent}

\newcommand{\leb}{\text{\rm Leb}}

\newcommand{\ser}{{\rm ser}}

\newcommand{\rank}{{\rm rank}}
\newcommand{\conv}{{\rm conv}}

\renewcommand{\proof}{\noindent{\bf Proof.\ }}

\begin{document}

\title[]{Rank-based stochastic differential inclusions and diffusion limits for a load balancing model}

\author{Rami Atar}
\address{Viterbi Faculty of Electrical and Computer Engineering,
Technion -- Israel Institute of Technology
} 
\author{Tomoyuki Ichiba}
\address{Department of Statistics and Applied Probability, University of California
Santa Barbara
}


\date{\today}

\begin{abstract}
In \cite{banerjee2023load}, a randomized load balancing model was studied in a heavy traffic asymptotic regime where the load balancing stream is thin compared to the total arrival stream. It was shown that the limit is given by a system of rank-based Brownian particles on the half-line. This paper extends the results of \cite{banerjee2023load} from the case of exponential service time to an invariance principle, where service times have finite second moment. The main tool is a new notion of rank-based stochastic differential inclusion, which may be of interest in its own right. 
\end{abstract}

\maketitle

\section{Introduction}

This paper studies a randomized load balancing model under an asymptotic regime introduced by Banerjee, Budhiraja and Estevez \cite{banerjee2023load}, where $N$ servers, labeled $1,\ldots,N$, cater to $N+1$ streams of jobs, labeled $0,\ldots,N$. For each $1\le i\le N$, server $i$ caters to stream $i$, and stream $0$ undergoes randomized load balancing. Specifically, the power-of-choice algorithm is applied to this stream, in which $\ell$ out of the $N$ queues are chosen at random and the job is routed to the shortest among them.  In the regime proposed in \cite{banerjee2023load}, the load balancing stream is much thinner than the remaining streams. The motivation to study a thin load balancing stream stems from the fact that the intensity of this stream, along with the parameter $\ell$, determines the communication volume between the dispatcher and servers, which one wants to keep low in practical applications.  In \cite{banerjee2023load}, $N$ is fixed, and, denoting the scaling parameter by $n$, the server capacity and the arrival rates of streams $1,\ldots,N$ scale like $n$, whereas the arrival rate of stream $0$ scales like $n^{1/2}$. Under a critical load condition, it is shown that the $N$ queue lengths, normalized by $n^{1/2}$, converge to a system of Brownian particles on the half line, with rank-dependent drift coefficients. This result thus identifies the minimal order of magnitude of the load balancing stream intensity at which load balancing has a macroscopic effect on the system behavior at the diffusion scale.

The treatment in \cite{banerjee2023load} assumes exponential service times.  The goal of this paper is to extend the diffusion limit result to an invariance principle, where service times are only assumed to possess a second moment. This does not amount to a mere technical improvement of the proof ideas of \cite{banerjee2023load}. To briefly explain this point, let $\hat X^n_i$, $1\le i\le N$ denote the normalized queue length process of server $i$. Let the {\it ranked} normalized queue lengths $\hat Y^n_i$, $1\le i\le N$, be defined as
\[
\hat Y^n_i(t)=\hat X^n_{\pi_t(i)}(t),
\]
for a $t$-dependent permutation $\pi_t$ of $\{1, \ldots , N\}$ that ensures that, for all $t$,
\[
\hat Y^n_1(t)\le\hat Y^n_2(t)\le\cdots\le\hat Y^n_N(t).
\]
The limiting dynamics of $\hat X^n_i$ is expected to be given in terms of a system of stochastic differential equations (SDEs) on $\R_+^N$ of the form (written here in a slightly simplified manner; see Section \ref{sec23} for a precise definition and the general form of SDEs treated here), 
\begin{equation}\label{103}
X_i(t)=X_i(0)+B_i(t)+\int_0^t b_{\calR_i(s)}ds+L_i(t),
\end{equation}
where $b_i$ are constants, $\calR_i(t)$ is the rank of $X_i(t)$ (e.g.,\ $\calR_i(t)=1$ if $X_i(t)=\min_jX_j(t)$), $B_i$ are mutually independent Brownian motions, and $L_i$ is a local time term for process $X_i$ at the origin. In particular, $b_{\calR_i(t)}$ are rank-dependent drift coefficients. If now one lets $Y_i$, $1\le i\le N$ denote the ranked processes corresponding to $\{X_i\}$  (by a transformation like the one above for the normalized queue lengths), one finds that they satisfy a system of SDEs of the form
\begin{equation}\label{104}
Y_i(t)=Y_i(0)+\tilde B_i(t)+b_it+\frac{1}{2}\tilde L_{i}(t)-1_{\{i<N\}}\frac{1}{2}\tilde L_{i+1}(t)+1_{\{i=1\}}\frac{1}{2}\tilde L_1(t),
\end{equation}
where $\tilde L_i$ are local time terms. For background on ranked processes of semimartingales, see \cite{MR2428716}. 
The strategy of the proof in \cite{banerjee2023load} is to work with
the approximating processes $\hat Y^n_i(t)$ and show convergence to the unique solution of \eqref{104}.  In this technique, one needs to show that the intersection time 
\[\int_{0}^{\cdot} 1_{\{ \hat Y_i^n(t)=\hat Y_{i+1}^n(t)\}} d t \] 
converges to zero in probability, as $n \to \infty$ for every $i = 1, \ldots , N-1$ (Lemma 4.2 and Corollary 4.3 of \cite{banerjee2023load}), and here the exponential service time assumption is convenient. Estimating it under more general conditions may be quite challenging.

Our proof strategy is to work with the processes $\hat X^n_i$ and show that they converge to the unique solution of \eqref{103}.
Because \eqref{103} itself does not involve multiple local-time terms except the one for spending time at zero, there is no need
to estimate the intersection times. On the other hand, our approach requires the notion of rank-based stochastic differential inclusions (SDIs), such as
\begin{equation}\label{105}
X_i(t)=X_i(0)+B_i(t)+\int_0^t\beta_i(s)ds+L_i(t).
\end{equation}
Here, $\beta_i(t)=b_{\calR_i(t)}$ at times when all particles are isolated from each other, but a milder condition is required at times when some of the particles collide. The precise condition is formulated by requiring $\{\beta_i(t)\}$ to belong to a set that depends on the current state $\{X_i(t)\}$
(see Section \ref{sec23}).

Differential inclusions play an important role in the study of non-smooth dynamical systems in both deterministic \cite{aub-cel, filippov2013differential} and stochastic \cite{atabudram} settings. An intuitive explanation for their effectiveness here is that while a precise formulation of \eqref{103} must specify the tie-breaking rule for how $\calR_i(t)$ are defined when two particles collide, such details should not matter. In fact, a differential inclusion  avoids keeping track of such details. When weak limits are taken for the queueing model, details on tie breaking rules are lost, and one is left with a differential inclusion which, in its simplest form, is given by \eqref{105}.

For a survey on randomized load balancing algorithms and their recent extensive study in asymptotic regimes, see \cite{der2022scalable}. Rank-based diffusions have also been extensively studied in recent years; see e.g., \cite{MR1894767, MR2187296, MR4752246, MR4420428, MR3771768, MR3672832, MR3055258}.

\subsection{Notation}
Throughout the paper, we use the following notation. $[N]=\{1,\ldots,N\}$, $\R_+=[0,\iy)$. $\iota$ denotes the identity map on $\R_+$. In $\R^N$, the Euclidean norm is denoted by $\|\cdot\|$. For $(X,d_X)$ a Polish space, let $C(\R_+,X)$ and $D(\R_+,X)$ denote the space of continuous paths and, respectively, \cadlag paths, endowed with the topology of uniform convergence on compacts and, respectively, the Skorokhod $J_1$ topology. Let $C_0^\up(\R_+,\R_+)$ denote the subset of $C(\R_+,\R_+)$ of nondecreasing functions that vanish at zero. For $f:\R_+\to\R^N$, denote
\begin{align*}
\|f\|^*_t&:=\sup_{s\in[0,t]}\|f(s)\|,\\
w_t(f,\del)&:=\sup\{\|f(u)-f(s)\|:u,s\in[0,t],|u-s|\le\del\}.
\end{align*}
If $v\in\R^N$, then $v_i$, $i\in[N]$ denote its coordinates, and vice versa: if $v_i$, $i\in[N]$ are given, then $v$ denotes $(v_1,\ldots,v_N)$. The same convention holds for random elements $v_i$ and stochastic processes $v_i(\cdot)$.
The symbol $\To$ denotes convergence in distribution. 

\section{Model and main results}\label{sec2}
\beginsec

\subsection{The load balancing model}
The model we consider is defined on a probability space $(\Om,\calF,\PP)$ with the corresponding expectation $\E$.  It consists of $N+1$ arrival streams, labeled $0,1,\ldots, N$, as well as $N$ queues and $N$ servers, labeled $1,\ldots,N$. For $i\in[N]$,
server $i$ serves jobs in queue $i$ in the order of arrival to the queue. The $N$ queues are fed by the $N+1$ arrival streams, where all jobs from arrival stream $i\in[N]$ are routed to queue $i$, whereas jobs from stream $0$, called the {\it load balancing stream} (LBS), are routed to the $N$ queues according to a randomized load balancing
algorithm described below.

The sequence of systems is indexed by $n\in\N$. The processes $X^n_i$, $E^n_i$, $D^n_i$ and $T^n_i$, $i\in[N]$, represent the $i$-th queue length process, arrival process, departure process and cumulative busyness process, respectively. Moreover, $A^n_0$ denotes the LBS arrival process and $A^n_i$ the process counting LBS arrivals routed by the algorithm to server $i$. 

For each $n \in \mathbb N$, the $N+1$ arrival processes $E^n_i$ and $A^n_0$ are mutually independent Poisson processes of intensities $\la^n_i$ and $\la^n_0$, respectively. Alternatively, they can be viewed as $N+1$ Poisson thinned streams of a single arrival stream, obtained by random selection.  These processes are assumed to have
right-continuous sample paths. We have the balance equation  
\begin{equation}\label{01}
X^n_i(t)=X^n_i(0-)+E^n_i(t)+A^n_i(t)-D^n_i(t),
\qquad i\in[N],\ t\in\R_+,
\end{equation}
and, assuming work conservation, the cumulative busyness process is given by
\begin{equation}\label{03}
T^n_i(t):=\int_0^t1_{\{X^n_i(s)>0\}}ds, \quad i \in [N], \, t \in \mathbb R_{+}. 
\end{equation}
Note that $\, t - T^{n}_{i}(t) \,$, $\, t \ge 0 \,$, $i \in [N]$ are the cumulative idle time processes.  

The load balancing algorithm routes LBS jobs to queues according to their relative lengths.
The precise construction requires the definition of a rank function. Namely,
$\rank: [N] \times \mathbb R^{N} \to [N]$ is defined as
\begin{equation}\label{17-}
\rank(i;x):=\#\{j \in [N] :x_j<x_i\}+\#\{j\le i: x_j=x_i\},
\qquad i\in[N],\,x\in\R^N. 
\end{equation}
Thus, $\rank(i;x)$ is the rank of $x_i$ among $\{x_j\}$, where the tie-breaking rule
is that the smaller index among the ties is preferred.
For example, for $\,x  = (1, 1, 2,2, 3) \in \mathbb R^{5}\,$, one has $\rank(i;x)_{i=1}^5=(1,2,3,4,5)$, and for $\,x = (1, 1, 3, 2, 2) \in \mathbb R^{5}\,$, $\rank(i;x)_{i=1}^5=(1,2,5,3,4)$.

Next, a probability vector $p=(p_r)\in[0,1]^N$, $p_1+\cdots+p_N=1$, is given,
assumed throughout to satisfy
\begin{equation}\label{010}
p_1\ge p_2\ge\cdots\ge p_N.
\end{equation}
The load balancing algorithm routes an LBS job to the queue whose length is ranked $r$ with probability $p_r$ (in particular, shorter queues are preferred). For the construction of this part of the queueing model, let $\theta_k$, $k\in\N$ be IID random variables with the common distribution $\PP(\theta_1=r)=p_r$, $r\in[N]$. 
Then the cumulative number $A^{n}_{i}(t)$ of the LBS arrivals routed to server $i$ by time $t$ is given by 
\begin{equation}\label{02}
A^n_i(t)=\int_{[0,t]}1_{\{\calR^n_i(s-)=\theta_{A^n_0(s)}\}}dA^n_0(s),
\qquad
\calR^n_i(t)=\rank(i;X^n(t)),
\qquad i\in[N],\ t\in\R_+.
\end{equation}

The most important special cases of this setting are two versions of the well-known {\it power-of-choice} algorithm. Here, upon an LBS arrival, the queue lengths of $\ell$ out of the $N$ queues, chosen uniformly at random, are sampled (with or without replacement).
The arrival is routed to the queue that is the shortest among the $\ell$ queues.
The volume of communication between the servers and the dispatcher is therefore
proportional to $\ell$. For this reason, in practice, $\ell$ is usually chosen much smaller than $N$.

Under the power-of-choice,
the probability $p_r$ that an LBS arrival is routed to the queue whose rank (as defined by \eqref{17-}) is $r$, is given by
\begin{equation} \label{eq:p*r}
p_r=\Big(\frac{N-r+1}{N}\Big)^\ell - \Big(\frac{N-r}{N}\Big)^\ell,
\qquad
\text{and}
\qquad
p_r=\frac{{N-r\choose \ell-1}}{{N\choose \ell}}, \quad \text{respectively}
\qquad r\in[N]
\end{equation}
for sampling with and without replacement, respectively.
Here, ${k\choose j}=0$ when $j>k$.
In both cases, $\sum_rp_{r}\, =\,  1$, and \eqref{010} hold. 
Our results are concerned with general $p$ satisfying \eqref{010},
but the main interest is in the cases \eqref{eq:p*r}.

\paragraph{\bf Initial conditions}
The residual times of jobs that have already been processed at time $0$ are assumed to satisfy some mild conditions. Denote the (random) set of queues that at time $0$ contains no jobs and, respectively, at least one job, by $\calN^n=\{i\in[N]:X^n_i(0)=0\}$ and $\calP^n=\{i\in[N]:X^n_i(0)>0\}$. Note that $\calN^{n}$ and $\calP^{n}$ partition $[N]$ with $\# \calN^{n} + \# \calP^{n} = N $. 

For $i\in\calP^n$, let $Z^n_i(0)$ denote the initial residual time of the head-of-line job in queue $i$. For $i\in\calN^n$, let fictitious jobs be added, having zero processing time. To this end, rather than specifying $X^n(0)$ as the initial queue length, $X^n(0-)$ is specified (as in \eqref{01}); and for each $i\in\calN^n$, the queue length is set to $X^n_i(0-)=1$ and the residual processing time is set to $Z^n_i(0)=0$. Note that this results in $X^n_i(0)=0$. Note that by adding the fictitious jobs we attain the following. The first job to enter service after time $0$ will have index $1$ regardless of the initial status of the queue (empty or non-empty). The initial condition is thus a tuple
\[
\calI^n=(\{X^n_i(0-),Z^n_i(0),i\in[N]\},\, \calN^n,\, \calP^n),
\]
where $(\calN^n,\calP^n)$ partitions $[N]$, and
\[
\begin{split}
X^n_i(0-)&=1,\ Z^n_i(0)=0,\ i\in\calN^n,
\\
X^n_i(0-)&\ge1,\ Z^n_i(0)>0,\ i\in\calP^n.
\end{split}
\]

\paragraph{\bf Service times}
Let $\Phi^\ser_i$, $i\in[N]$, be Borel probability measures on $[0,\iy)$ with mean 1, standard deviation $\sig^\ser_i\in(0,\iy)$, and $\Phi^\ser_i(\{0\})=0$. Let $\Phi^{\ser,n}_i$ be defined as scaled versions of these measures, uniquely specified via
\begin{equation}
\Phi^{\ser,n}_i[0,x]:=\Phi^\ser_i[0,\mu^n_ix], \quad x\in\R_+.
\end{equation}
Here, $\mu^n_i>0$ is the service rate of server $i$ in the $n$-th system.
For $k\ge1$, let $Z^n_i(k)$ denote the service time
of the $k$-th job to be served by server $i$ after the head-of-line job at time $0-$ (for $i\in\calN^n$ this means the $k$-th job after the fictitious one). It is assumed that, for every $i$, $\calZ^n_i=(Z^n_i(k),k\ge1)$ is an IID sequence with common distribution $\Phi^{\ser,n}_i$.

The potential service process $S^n_i$, evaluated at $t$, represents the number of jobs completed by server $i$ by the time it has worked $t$ units of time. With $\sum_0^{-1}=0$, it is given by 
\begin{equation}\label{34}
S^n_i(t)=
\max\Big\{k\in\Z_+:\sum_{j=0}^{k-1}Z^n_i(j)\le t\Big\}, \qquad t\ge0.
\end{equation}
The departure processes are, therefore, given by $D^n_i(t)=S^n_i(T^n_i(t))$. This is the number of jobs completed by time $t$ by server $i$.

It is assumed that, for each $n$, the $2N+3$ stochastic elements
\[
E^n_i, i\in[N], \quad \calZ^n_i, i\in[N], \quad \calI^n,
\quad A^n_0, \quad \{\theta_k\},
\]
are mutually independent.

\subsection{\bf The scaling and critical load condition.}
The arrival and service rates are assumed to satisfy the following. There are constants $\la_i>0$ and $\hat\la_i\in\R$ such that
\begin{equation}\label{04}
\hat\la^n_i:=n^{-1/2}(\la^n_i-n\la_i)\to\hat\la_i , \quad \text{ as } n\to\iy,
\qquad i\in[N],
\end{equation}
a constant $\la_0>0$ such that
\begin{equation}\label{05}
\hat\la^n_0:=n^{-1/2}\la^n_0\to \la_0 , \quad \text{ as } n\to\iy,
\end{equation}
and constants $\mu_i>0$ and $\hat\mu_i\in\R$ such that
\begin{equation}\label{06}
\hat\mu^n_i:=n^{-1/2}(\mu^n_i-n\mu_i)\to\hat\mu_i , \quad \text{ as } n\to\iy,
\qquad i\in[N].
\end{equation}
Each of the queues is critically loaded, namely
\begin{equation}\label{09}
\la_i\equiv\mu_i, \qquad i\in[N].
\end{equation}
The scaled initial residual time $n^{1/2} Z^{n}_{i}(0) $
are assumed to satisfy
\begin{equation}\label{40}
n^{1/2}Z^n_i(0)\to0 \quad \text{ in probability,} \qquad i\in[N], 
\end{equation}
as $n \to \infty$.

Finally,  we specify conditions regarding the initial queue lengths $X^n_i(0-)$, representing two different scenarios.  One with initial queue lengths of the order $n^{1/2}$, and another with queue lengths of a larger order of magnitude. To state these conditions, let us define
\[
\hat X^n(0-) = (\hat X^n_{1}(0-),\ldots,\hat X^n_{N}(0-)),
\quad
\hat X^n_{i}(0-) := n^{-1/2} X_{i}^{n}(0-) ,\, \,  i \in [N],
\]
and, for a fixed sequence $\al_n$ satisfying $n^{-1/2}\al_n\to\iy$,
\[
\check X^n(0-) = (\check X^n_{1}(0-),\ldots,\check X^n_{N}(0-)),
\quad
\check X^n_{i}(0-) := n^{-1/2} (X_{i}^{n}(0-)-\al_n) ,\, \,  i \in [N].
\]
It will be assumed that one of the following holds: Either
\begin{equation}\label{IC0}\tag{IC${}_0$}
\hat X^n(0-) \,  \To \, X_{0} := (X_{0, 1}, \ldots , X_{0, N}) ,
\text{ an $\R_{+}^{[N]}$-valued random vector,}
\end{equation}
or
\begin{equation}\label{ICalpha}\tag{IC${}_\al$}
\check X^n(0-) \,  \To \, X_{0} := (X_{0, 1}, \ldots , X_{0, N}) ,
\text{ an $\R^{[N]}$-valued random vector.}
\end{equation}

The main results of this model are concerned with diffusion-scale versions of the queue length processes $X^n_i(t)$, $t \ge 0 $ and cumulative idle time processes $t - T_{i}^{n}(t)$, $t \ge 0 $.
In particular, under \eqref{IC0}, we will be interested in
\begin{equation}\label{07}
\hat X^n_i(t):=n^{-1/2}X^n_i(t),
\qquad
\hat L^n_i(t):=n^{-1/2}\mu^n_i(t-T^n_i(t)), \quad t \ge 0 , \, \, i \in [N],
\end{equation}
whereas under \eqref{ICalpha}, we will study the behavior of
\begin{equation}\label{07a}
\check X^n_i(t):=n^{-1/2}(X^n_i(t)-\al_n), \quad t \ge 0 , \, \, i \in [N],
\end{equation}
where the queue lengths are centered around the same constants
$\al_n$ as the initial conditions.

\subsection{Rank-based SDE and SDI}\label{sec23}
We are concerned with a rank-based diffusion defined via
a system of SDEs or SDIs, with and without reflection.

Let constants $b := (b_{1}, \ldots, b_{N} ) \in\R^N$, $m := (m_{1}, \ldots , m_{N}) \in\R^N$, $\sig := (\sigma_{1}, \ldots , \sigma_{N}) \in(0,\iy)^N$ be given, and consider
rank-based SDE without, and, respectively, with reflection,
\begin{equation}\label{SDE}\tag{SDE}
\begin{split}
&X_i(t)=X_{0,i}+\sig_iB_i(t)+m_it+\int_0^tb_{\calR_i(s)}ds,
\\
&\calR_i(t)=\rank(i;X(t)), 
\end{split}
\hspace{6em} t \ge 0 , \, \, i \in [N] .
\end{equation}
\begin{equation}\label{SDER}\tag{SDER}
\begin{split}
&X_i(t)=X_{0,i}+\sig_iB_i(t)+m_it+\int_0^tb_{\calR_i(s)}ds+L_i(t)\ge0,
\\
&\calR_i(t)=\rank(i;X(t)), 
\\
&\int_0^\iy X_i(t)dL_i(t)=0,
\end{split}
\hspace{2em} t \ge 0 , \, \, i \in [N] .
\end{equation}
Here, $\rank$ is the function defined in \eqref{17-} and $L_{i}$ are the continuous, nondecreasing, adapted processes, starting from $L_{i}(0) = 0$, that make the $N$-dimensional process $X(\cdot)$ stay in the non-negative orthant $\mathbb R_{+}^{N}$. 

The precise notions of a solution, in strong and weak form, are standard, but we give them for completeness.
On a given stochastic basis $(\Omega, \mathcal F, \mathbb F := \{\mathcal F_{t}, t \ge 0 \}, \mathbb P) $ satisfying the usual conditions,
with an $\mathbb F$-Brownian motion $B$ in dimension $N$, and an initial condition $X_0\in\calF_0$, a strong solution of \eqref{SDER} is an $\mathbb F$-adapted process $(X, L)$ with sample paths in $C(\R_+,\R_+^N)\times C_0^\up(\R_+,\R_+)^N$ that satisfies \eqref{SDER}. A weak solution to \eqref{SDER} is 
a stochastic basis $(\Omega, \mathcal F, \mathbb F, \mathbb P) $ that satisfies the usual conditions, along with processes $(X,L,B)$ defined on it, satisfying \eqref{SDER} a.s., where $B$ is an $\mathbb F$-Brownian motion in dimension $N$,
$X_0\in F_0$, and $(X,L)$ is $\mathbb F$-adapted, with sample paths in $C(\R_+,\R_+^N)\times C_0^\up(\R_+,\R_+)^N$.

With a slight abuse of terminology, we will sometimes say that a tuple $(X,L,B)$,
or even $(X,L)$ is a weak solution, without specifying the stochastic basis
(or the Brownian motion).

We say that uniqueness in distribution holds for \eqref{SDER} if for any two weak solutions $(X,L,B)$,  $(\Omega, \mathcal F, \mathbb F, \mathbb P) $ and $( \widetilde{X},\widetilde{L},\widetilde{B})$,  $(\widetilde{ \Omega}, \widetilde{ \mathcal F}, \widetilde{ \mathbb F}, \widetilde{ \mathbb P}) $, with the same initial distribution, the two processes $(X,L)$ and $ (\widetilde{X}, \widetilde{L})$ have the same distribution. We say that pathwise uniqueness holds for \eqref{SDER}, if for any two weak solutions $(X,L,B)$,  $(\Omega, \mathcal F, \mathbb F, \mathbb P) $ and $( \widetilde{X},\widetilde{L},{B})$,  $({ \Omega}, { \mathcal F}, { \mathbb F} , { \mathbb P}) $,
with common initial value, the two processes $X$ and $ \widetilde{X}$ are indistinguishable, that is, $\mathbb P (X_{t} = \widetilde{X}_{t}, t \ge 0 ) \, =\,  1 $,
and so are $L$ and $\widetilde{L}$.

Analogous notions are defined for \eqref{SDE} in a similar manner,
with $X$ having the sample path in $C(\R_+,\R^N)$. 

\begin{remark} \label{rem: Yamada-Watanabe-Kurtz}
    It is shown by Yamada and Watanabe in \cite{MR278420} that the pathwise uniqueness and the weak existence of SDE imply the strong existence of the solution and the uniqueness in law. 
    Following the weak and strong solutions of stochastic system developed by Kurtz in \cite{MR3254737}, we claim that the pathwise uniqueness and the weak existence of SDI imply the strong existence of the solution and the uniqueness in law. 
Particularly, the pathwise uniqueness and weak existence of \eqref{SDIR} $($\eqref{SDI}, respectively$)$ imply the strong existence of the solution of \eqref{SDIR} $($\eqref{SDI}, respectively$)$ and uniqueness in law.  
\end{remark}

To introduce our notion of a rank-based SDI,
let $\Pi$ denote the set of all permutations of $[N]$.
For each $\pi\in\Pi$, write $b_\pi$ for the vector $(b_{\pi(1)},\ldots,b_{\pi(N)})$. 
Let us consider a set-valued map $\mathfrak P:\mathbb R^{N} \to 2^{\Pi}$ defined by 
\begin{equation} \label{eq:Px}
\mathfrak P(x)=\{\pi\in\Pi:x_i<x_j \text{ implies } \pi(i)<\pi(j) \text{ for every } i, j \in [N] 
\}, 
\end{equation}
for $x \in \R^{N}$.
Denote by $\text{conv}(A)$ the convex hull of set $A\subset\R^N$. 
The SDI without and, respectively, with reflection, that will be of interest here, are
\begin{equation}\label{SDI}\tag{SDI}
\begin{split}
&X_i(t)=X_{0,i}+\sig_iB_i(t)+m_it+\int_0^t\beta_i(s)ds,\hspace{5.4em} t\ge0,\ i\in[N],
\\
&\beta(t)=(\beta_1(t),\ldots,\beta_N(t))\in\conv\{b_\pi:\pi\in{\mathfrak P}(X(t))\} \qquad \text{a.e. } t \in \mathbb R_{+}, 
\end{split}
\end{equation}
and a similar one on the state space $\R_+^N$,
\begin{equation}\label{SDIR}\tag{SDIR}
\begin{split}
&X_i(t)=X_{0,i}+\sig_iB_i(t)+m_it+\int_0^t\beta_i(s)ds+L_i(t),\hspace{2em} t\ge0,\ i\in[N],
\\
&\beta(t)\in\conv\{b_\pi:\pi\in\mathfrak P(X(t))\} \hspace{10.7em} \text{a.e. } t \in \mathbb R_{+} ,
\\
&\int_{0}^{\infty} X_i(t) dL_i(t) = 0 , \hspace{14.5em} i \in [N] .
\end{split}
\end{equation}

We call $\beta$ the {\it rank-dependent drift process}.
We say that $X$ (respectively, $(X,L)$) is a solution to \eqref{SDI}
(respectively, \eqref{SDIR}) if there exist an $\mathbb F$-progressively measurable process $\beta$ so that \eqref{SDI}
(respectively, \eqref{SDIR}) holds.
We extend the notions of weak/strong solution, weak uniqueness and pathwise uniqueness, given above, to \eqref{SDI} and \eqref{SDIR}.

\subsection{Results}

A vector $b\in\R^N$ is said to be nonincreasing if the sequence $b_1,\ldots,b_N$ is.

\begin{theorem}\label{th1}
Let data $b,m\in\R^N$, $\sig\in(0,\iy)^N$ be given and assume that $b$ is nonincreasing. Then pathwise uniqueness and strong existence hold
for \eqref{SDE}, \eqref{SDER}, \eqref{SDI} and \eqref{SDIR}.
\end{theorem}

\begin{theorem}\label{th2}
Consider data defined in terms of the load balancing model, as follows
\begin{equation}\label{011}
b_r=\la_0p_r, \ r\in[N],
\qquad
m_i=\hat\la_i-\hat\mu_i,
\qquad
\sig_i=(\la_i+\mu_i(\sig_i^\ser)^2)^{1/2}, \qquad i\in[N].
\end{equation}
i. Assume \eqref{IC0}. Then $(\hat X^n,\hat L^n)\To(X,L)$ in $D(\R_+,\R^N)\times C(\R_+,\R^N)$ as $n\to\iy$, where $(X,L)$ is the solution to \eqref{SDIR}, equivalently, \eqref{SDER}, with data \eqref{011}.
\\
ii. Fix a sequence $\al_n$ with $n^{-1/2}\al_n\to\iy$ and assume \eqref{ICalpha}. Then $(\check X^n,\hat L^n)\To(X,0)$ in $D(\R_+,\R^N)\times C(\R_+,\R^N)$ as $n\to\iy$, where $X$  is the solution to \eqref{SDI}, equivalently, \eqref{SDE}, with data \eqref{011}. 
\end{theorem}

\subsection{Proof outline}\label{sec25}
The proofs of Theorems \ref{th1} and \ref{th2} are intertwined:
The existence of a limit for the queueing model is based on the uniqueness
provided by Theorem \ref{th1}, whereas existence of solutions to both differential inclusions is a consequence of the convergence
proved in Theorem \ref{th2}. The steps are as follows.

\noi
1. Pathwise uniqueness of \eqref{SDI} and \eqref{SDIR}.
This is shown in Section \ref{sec3}, specifically, in Proposition \ref{prop0}.

\noi
2. This automatically gives the pathwise uniqueness of \eqref{SDE}
and \eqref{SDER}. See Remark \ref{rem0}.

\noi
3. Weak convergence of the rescaled queueing model to a solution of \eqref{SDIR} and \eqref{SDI} (under \eqref{IC0} and \eqref{ICalpha}, respectively). This is argued by showing that tightness holds and that subsequential weak limits satisfy the differential inclusions, which, in view of step 1, imply the existence of a limit of the entire sequence. This is carried out in Section \ref{sec4}.

\noi
4. Step 3 immediately gives weak existence of solutions to
both differential inclusions.

\noi
5. For \eqref{SDI} and \eqref{SDIR}, the set of times when two or more components $X_i$ meet
is shown to have Lesbegue measure zero. This gives weak existence of solutions to \eqref{SDE} and \eqref{SDER}.
This is proved in Section \ref{sec5}.

\noi
6. The Yamada-Watanabe Theorem now gives strong existence for the four equations \eqref{SDI}, \eqref{SDIR}, \eqref{SDE} and \eqref{SDER} (cf. Remark \ref{rem: Yamada-Watanabe-Kurtz}). This is also proved in Section \ref{sec5}.

\section{SDI and SDE uniqueness}\label{sec3}
\beginsec

\begin{proposition}\label{prop0}
Pathwise uniqueness holds for \eqref{SDI} and \eqref{SDIR}.
\end{proposition}
\begin{remark}\label{rem0}
Note that every solution to \eqref{SDE} is a solution to \eqref{SDI}. The same holds for \eqref{SDER} and \eqref{SDIR}. Hence,  the above immediately gives the pathwise uniqueness
to \eqref{SDE} and \eqref{SDER}.
\end{remark}

\noi{\bf Proof.}
The starting point for this proof is the {\it rearrangement inequality}
\cite{hardy1952inequalities}, which states that if $u_i$, $i\in[N]$ is nonincreasing and $v_i$, $i\in[N]$ is nondecreasing, then 
\begin{equation}\label{RearrIneq}
\sum_{i=1}^Nu_iv_i\le\sum_{i=1}^Nu_{\pi(i)}v_i 
\end{equation}
for any permutation $\pi \in \Pi$. 

Let two weak solutions $(X,L,B)$ and $(Y,M,B)$, defined on the same
stochastic basis, be given, and let $\beta$ and $\gamma$ be the corresponding
rank-dependent drift processes.
The difference process $V=(V_{1}, \ldots , V_{N}) $ with $V_{i}(\cdot) := X_{i}(\cdot) - Y_{i}(\cdot) $, $i \in [N]$ satisfies 
\[
V_i(\cdot):=X_i(\cdot) -Y_i(\cdot)=\int_0^\cdot(\beta_i-\gam_i)ds
+L_i (\cdot) -M_i (\cdot), 
\]
and hence, the squared norm satisfies
\begin{align*}
\frac{1}{2}\|V(t)\|^2&=\int_0^tV(s)\cdot[(\beta(s)-\gam(s))ds+dL(s)-dM(s)]
\\
&=\sum_{i=1}^{N}\int_0^t(X_i(s)-Y_i(s))(\beta_i(s)-\gam_i(s))ds
+\sum_{i=1}^{N}\int_0^t(X_i(s)-Y_i(s))(dL_i(s)-dM_i(s)).
\end{align*}
The proof of pathwise uniqueness for \eqref{SDIR} will be complete if one shows that the right-hand above is non-positive for all $t$.
 The second sum is clearly nonpositive because $\int_{0}^{\cdot} X_i(t) dL_i(t) =0$ while $\int_{0}^{\cdot} X_i(t) dM_i(t) \ge0$, and similarly for the $Y_i$ parts in the second sum. Thus, it suffices to show that
 \begin{equation} \label{eq:XYbg}
\sum_{i=1}^{N} ( X_{i} (t) - Y_{i}(t) ) ( \beta_{i} (t) - \gamma_{i} (t))\le0 , \quad \text{a.e. } \, t \in \mathbb R_{+} .
 \end{equation}
A similar argument for \eqref{SDI}, only slightly simpler as it does not involve
the boundary terms $L$ and $M$, also leads to the conclusion that
proving \eqref{eq:XYbg} suffices. We therefore proceed to show \eqref{eq:XYbg}
for both \eqref{SDI} and \eqref{SDIR}.

Following the definition of the SDI, let us write 
\[
\beta(t)=\sum_{\pi \in \Pi} g_\pi(t) b_\pi , \quad \gamma(t)=\sum_{\pi \in \Pi} h_\pi(t) b_\pi , \quad t \ge  0 , 
\]
where 
$g_\pi(t)\ge0$, $\pi \in \Pi$, $\sum_{\pi \in \Pi} g_\pi(t)=1$, and $g_\pi(t)=0$ for $\pi\notin\mathfrak P(X(t))$ and similarly, 
$h_\pi(t)\ge0$, $\pi \in \Pi$, $\sum_{\pi \in \Pi} h_\pi(t)=1$, and $h_\pi(t)=0$ for $\pi\notin\mathfrak P(Y(t))$ 
for $t \ge 0 $. 
By the definition of $\mathfrak P (\cdot)$ in \eqref{eq:Px}, we have therefore that for a.e. $\,t$, the conjunction of conditions  $g_\pi(t)>0$ and $X_i(t)<X_j(t)$ implies $\pi(i)<\pi(j)$. A similar statement holds for $h_{\pi}(t)$ and $Y(t)$. Then, omitting $t$ in the formulas below, we now have 
\begin{align*}
\sum_{i=1}^{N}(X_i-Y_i)(\beta_i-\gam_i)
&=\sum_{i=1}^{N}(X_i-Y_i)\sum_{\pi \in \Pi} [g_\pi b_{\pi(i)}-h_\pi b_{\pi(i)}]
\\
&=\sum_{i=1}^{N}(X_i-Y_i)\sum_{\pi,\sig \in \Pi}g_\pi h_\sig(b_{\pi(i)}-b_{\sig(i)}) \\
&= \sum_{\pi,\sig \in \Pi} \sum_{i=1}^{N}(X_i-Y_i)g_\pi h_\sig(b_{\pi(i)}-b_{\sig(i)}).
\end{align*}
Consider the terms involving only the $X_i$'s in the above sum, namely
\begin{align*}
\sum_{\pi,\sig}\sum_{i=1}^{N}X_i g_\pi h_\sig(b_{\pi(i)}-b_{\sig(i)}).
\end{align*}
To prove nonpositivity of this sum, it suffices to prove that,
for any $\pi$ for which $g_\pi>0$
(that is, $\pi\in\mathfrak P(X)$), and for an arbitrary $\sig\in\Pi$,
\[
\sum_{i=1}^{N}X_i(b_{\pi(i)}-b_{\sig(i)})\le0.
\]
Because of the monotonicity: if $\pi(i)>\pi(j)$ then $X_i\ge X_j$, 
reordering $i$'s so that $\pi(i)$ is increasing will give that $X_i$ are nondecreasing.
It will also give (by the assumption on $b$) that $b_{\pi(i)}$ is nonincreasing. 
Hence, the last display follows from the rearrangement inequality \eqref{RearrIneq}. Interchanging the roles of $X,\pi, g_{\pi}$ and those of $Y, \sigma, h_{\sigma}$, we also obtain the inequality for $Y$.  Therefore, we conclude
that \eqref{eq:XYbg} is true.
This shows that $\,\lVert V(\cdot) \rVert \equiv 0 \,$ for both \eqref{SDI} and \eqref{SDIR}. 
\qed

\section{Diffusion limits}\label{sec4}
\beginsec

This section provides the main step toward proving
Theorem \ref{th2}. In Subsection \ref{sec41},
under condition \eqref{IC0}, it is shown that $(\hat X^n,\hat L^n)$
converge to the unique solution of \eqref{SDIR}.
This proves Theorem \ref{th2}(i) except for the statement regarding
\eqref{SDER}, which is treated later in Section \ref{sec5}.
Subsection \ref{sec42} assumes \eqref{ICalpha},
and, similarly, proves Theorem \ref{th2}(ii) except
the statement on \eqref{SDE}, whose proof is
deferred to Section \ref{sec5}.
Most of the work is done in Subsection \ref{sec41}.

\subsection{The limit under condition \texorpdfstring{\eqref{IC0}}{(IC0)}}
\label{sec41}

In this subsection we prove Proposition \ref{prop1} below. The proof relies on the uniqueness of \eqref{SDIR}
proved in the previous section. Let
\begin{equation}\label{07+}
\hat E^n_i(t)=n^{-1/2}(E^n_i(t)-\la^n_it),
\qquad
\hat S^n_i(t)=n^{-1/2}(S^n_i(t)-\mu^n_it),
\end{equation}
\begin{equation}\label{08}
\hat A^n_0(t)=n^{-1/2}A^n_0(t),
\qquad
\hat A^n_i(t)=n^{-1/2}A^n_i(t),
\end{equation}
\begin{equation}\label{39}
\hat P^{\#,n}_i(t)=\la_0\int_0^tp_{\calR^n_i(s)}ds,
\end{equation}
\begin{equation}\label{35}
\hat m^n_i=\hat\la^n_i-\hat\mu^n_i.
\end{equation}
Throughout this subsection, let condition \eqref{IC0} hold.

\begin{proposition}\label{prop1}
i.
The sequence $(\hat X^n,\hat L^n,\hat E^n,\hat S^n,\hat P^{\#,n})$ is $C$-tight.
\\
ii. If $(X,L,E,S,P)$ is a subsequential weak limit,
then $(X,L,B)$ forms a weak solution to \eqref{SDIR}
with the data indicated in Theorem \ref{th2}, and
where $B_i=\sig_i^{-1}(E_i-S_i)$ (with $\sig_i$ as in \eqref{011})
and the progressively measurable
rank-dependent drift $\beta$ is the a.e.\ derivative of $P$.
\\
iii. Consequently, denoting $\hat B^n_i=\sig_i^{-1}(\hat E^n_i-\hat S^n_i)$,
one has $(\hat X^n,\hat L^n,\hat B^n)\To(X,L,B)$, where the latter
is a (weak) solution of \eqref{SDIR}.
\end{proposition}

Let the Skorokhod map on the half-line
$\Gam: D(\R_+,\R)\to D(\R_+,\R_+)^2$ be defined by
\[
\Gam(y)=(x,z) \qquad \text{where}\qquad
x(t)=y(t)+z(t),
\qquad
z(t)=\sup_{s\in[0,t]}y^-(s),
\qquad
t\ge0.
\]
Note that if $(x,z)=\Gam(y)$ then
\begin{equation}\label{012}
z(t)\le\|y\|^*_t,
\qquad
w_t(z,\del)\le w_t(y,\del),
\qquad
t>0,\,\del>0.
\end{equation}

\begin{lemma}\label{lem:sk}
Let $(y,x,z)\in D(\R_+,\R)\times D(\R_+,\R_+)^2$ satisfy $x=y+z$.
Then the condition $z$ is nondecreasing and
$\int_{[0,\iy)}x(t)dz(t)=0$ (with the convention $z(0-)=0$)
holds if and only if $(x,z)=\Gam(y)$.
\end{lemma}
\proof
This is known as Skorokhod's lemma \cite[\S 8]{chu-wil}.
\qed

The tuple
\[
\calS^n(t)=(E^n_i(t),A^n_i(t),D^n_i(t),X^n_i(t),T^n_i(t),
\,i\in[n],\, A^n_0(t),\theta_{A^n_0(t)})
\]
is used to define the `history' of the system, namely the filtration
\[
\calF^n_t=\sig\{\calI^n,\calS^n(s),s\in[0,t]\}, \quad t \ge 0 . 
\]
Let
\begin{equation}\label{36}
\hat P^n_i(t)=\hat\la^n_0\int_0^tp_{\calR^n_i(s)}ds
\qquad \text{and}\qquad
\hat M^n_i(t)=\hat A^n_i(t)-\hat P^n_i(t).
\end{equation}
\begin{lemma}\label{lem1}
The process $\hat M^n_i$ is an $\{\calF^n_t\}$-martingale, with optional quadratic variation
$[\hat M^n_i](t)=n^{-1}A^n_i(t)$.
\end{lemma}

\proof
Clearly $A^n_i$ and $\calR^n_i$, defined in \eqref{02}, are $\calF^n_t$-adapted,
and $A^n_i(t)$ is integrable for all $t$. Hence, the same is true for $\hat M^n_i$.
Next, let
\[
t^n(k)=\inf\{t\ge0:A^n_0(t)\ge k\},\qquad k=1,2,\ldots.
\]
These are the stopping times on $\{\calF^n_t\}$. Hence
\begin{equation}\label{38}
t^n(k)\in\calF^n_{t^n(k)-}, \qquad k\ge1,
\end{equation}
where we recall that for a stopping time $\tau$,
\[
\calF^n_{\tau-}=\calF^n_0\vee\sig\{A : A\cap\{\tau<t\} \in \mathcal F^{n}_{t} , t \ge 0 \}
\]
(see \cite[I.1.11 and I.1.14]{jacshi}).
To show the martingale property, we can write, using
$\int_0^tp_{\calR^n_i(s-)}ds=\int_0^tp_{\calR^n_i(s)}ds$,
\begin{align*}
n^{1/2}\hat M^n_i(t)&=\int_{[0,t]}1_{\{\calR^n_i(s-)=\theta^n_{A^n_0(s)}\}}dA^n_0(s)
-\la^n_0\int_0^tp_{\calR^n_i(s)}ds
\\
&=\int_{[0,t]}(1_{\{\calR^n_i(s-)=\theta^n_{A^n_0(s)}\}}
-p_{\calR^n_i(s-)})dA^n_0(s)
+\int_0^tp_{\calR^n_i(s-)}(dA^n_0(s)-\la^n_0ds)
\\
&=: M^n_{i,1}(t)+M^n_{i,2}(t).
\end{align*}
For $M^n_{i,1}$, write
\[
A^n_i(t)=\sum_{k=1}^{A^n_0(t)}
1_{\{\calR^n_i(t^n(k)-)=\theta^n_k\}}.
\]
The history of the system up to $t^n(k)-$, namely $\{\calS^n(t),t<t^n(k)\}$, can be recovered from the tuple $\calI^n$, $(E^n_i(t),t\in\R_+,i\in[N])$, $(A^n_0(t), t\in\R_+)$, $(Z^n_i(j), j\in\N,i\in[N])$ and finally, $(\theta^n_j, j\le k-1)$. By our assumptions, $\theta^n_k$ is independent of this tuple. As a result, it is independent of $\calF^n_{t^n_i(k)-}$.
Therefore, for $0\le s<t$, we have
\begin{align*}
\E[A^n_i(t)|\calF^n_s]-A^n_i(s)
&=\sum_{k=1}^\iy\E[1_{\{\calR^n_i(t^n(k)-)=\theta^n_k\}}
1_{\{s<t^n(k)\le t\}}|\calF^n_s]
\\
&=
\sum_{k=1}^\iy\E[\E[1_{\{\calR^n_i(t^n(k)-)=\theta^n_k\}}
1_{\{s<t^n(k)\le t\}}|\calF^n_{t^n(k)-}]|\calF^n_s]
\\
&=
\sum_{k=1}^\iy\E[p_{\calR^n_i(t^n(k)-)}
1_{\{s<t^n(k)\le t\}}|\calF^n_s]
\\
&=
\E[C^n_i(t)|\calF^n_s]-C^n_i(s),
\end{align*}
where
\[
C^n_i(t)=\sum_{k=1}^{A^n_0(t)}p_{\calR^n_i(t^n(k)-)}
=\int_{[0,t]}p_{\calR^n_i(s-)} dA^n_0(s),
\]
showing that $A^n_i-C^n_i=M^n_{i,1}$ is a martingale.

In the expression for $M^n_{i,2}$, the integrand is $\{\calF^n_t\}$-adapted and has LCRL sample paths, while the integrator is a martingale on this filtration. As a result, $M^n_{i,2}$ is a local martingale \cite[Theorem II.20]{protter}; Using  the estimate $\|M^n_{i,2}\|^*_t\le A^n_0(t)+c$ shows that it is, in fact, a martingale. As a result, so is $\hat M^n_i$. Finally, the expression for the quadratic variation is straightforward.
\qed

\begin{lemma}\label{lem4}
i. One has
\begin{align}\label{10}
\hat X^n_i&=\hat U^n_i+\hat L^n_i \qquad \text{where}\qquad
\hat U^n_i(t)=\hat X^n_i(0-)+\hat E^n_i(t)+\hat A^n_i(t)
-\hat S^n_i(T^n_i(t))+\hat m^n_it.
\end{align}
Moreover, the sample paths of $\hat L^n_i$ are in $C^\up_0$ and
\begin{equation}\label{11}
\int_0^\iy\hat X^n_i(t)d\hat L^n_i(t)=0.
\end{equation}
In particular,
\begin{equation}\label{11+}
(\hat X^n_i,\hat L^n_i)=\Gam(\hat U^n_i), \qquad i\in[N].
\end{equation}

ii.
One has $(\hat E^n,\hat S^n)\To(E,S)$, where the latter is a pair of mutually independent $N$-dimensional Brownian motions starting at zero, with zero drift and diffusion coefficients $\diag(\la_i^{1/2})$ and $\diag(\mu_i^{1/2}\sig_i^\ser)$, respectively (where we recall $\la_i=\mu_i$).

iii.
$\hat P^n_i$, $\hat L^n_i$ and $\hat X^n_i$ are $C$-tight,
and $\hat M^n_i\to0$ in probability.
\end{lemma}

\proof
i. By \eqref{01} and \eqref{07},
\begin{align*}
n^{-1/2}X^n_i(t)&=n^{-1/2}X^n_i(0-)+n^{-1/2}(E^n_i(t)-\la^n_it)
+n^{-1/2}(\la^n_i-n\la_i)t+n^{1/2}\la_i t+n^{-1/2}A^n_i(t)\\
&\quad
-n^{-1/2}(S^n_i(T^n_i(t))-\mu^n_iT^n_i(t))-n^{-1/2}\mu^n_iT^n_i(t),
\end{align*}
and
\[
-n^{-1/2}\mu^n_iT^n_i(t)=-n^{-1/2}(\mu^n_i-n\mu_i)t-n^{1/2}\mu_i t+\hat L^n_i(t).
\]
Using \eqref{09}, \eqref{07+}, \eqref{08} and \eqref{35} gives \eqref{10}. The properties of $\hat L^n_i$ and \eqref{11} follow from \eqref{03}. The identity \eqref{11+} follows from Lemma \ref{lem:sk}.

ii.
The fact that $\hat E^n\To E$ follows from the central limit theorem for renewal processes \cite[\S 17]{bil}, and the fact that, by \eqref{04} $n^{-1}\la_i^n\to\la_i$, for each $i$. For $\hat S^n_i$, one has to be careful about the fact that the assumptions about $Z^n_i(0)$ differ from those about $Z^n_i(k)$, $k\ge1$. By \eqref{34}, $S^n_i$ is the inverse of
\[
Z^n_i(0)+\sum_{j=1}^{k-1}Z^n_i(j).
\]
In the expression
\[
n^{1/2}(Z^n_i(0)-(\mu^n_i)^{-1})
+n^{1/2}\sum_{j=1}^{[(n-1)t]}(Z^n_i(j)-(\mu^n_i)^{-1}),
\]
the first term converges to $0$ in probability by \eqref{40}.
In the second term, the summands are IID, and thus its limit in law
is a zero drift Brownian motion with diffusion $\sig^\ser_i$.
Hence again by  \cite[\S 17]{bil},
\eqref{06} and the independence of $\{Z^n_i(j)\}$ across $i$, we have $\hat S^n\To S$.
The mutual independence of $E$ and $S$ follows from that of
$\hat E^n$ and $\hat S^n$.

iii.
Because $\hat\la^n_0\to\la_0$, the processes $\hat P^n_i$
are all $(\la_0+1)$-Lipschitz, null at zero, for $n$ sufficiently large.
Hence, they are $C$-tight.
Moreover, by Lemma \ref{lem1} and the calculation
\[
\EE[n^{-1}A^n_i(t)]=n^{-1}\la^n_0\int_0^t\EE[p_{\calR^n_i(s)}]ds
\le n^{-1}(\la_0+1)n^{1/2}t,
\]
one has $\hat M^n_i\to0$ in probability. By \eqref{36}, this shows that $\hat A^n_i$ are also $C$-tight.

In view of the identity $(\hat X^n_i,\hat L^n_i)=\Gam(\hat U^n_i)$ and \eqref{012},
\begin{equation}\label{16}
\hat L^n_i(t)\le\|\hat U^n_i\|^*_t,\qquad w_t(\hat L^n_i,\del)\le w_t(\hat U^n_i,\del).
\end{equation}
Because $T^n_i$ are $1$-Lipschitz and $\hat S^n_i$ are $C$-tight,
so are $\hat S^n_i(T^n_i(\cdot))$. We have already shown that
$\hat E^n_i$ and $\hat A^n_i$ are $C$-tight.
Hence, by the convergence in law of $\hat X^n_i(0-)$
assumed in \eqref{IC0} and the convergence $\hat m^n_i\to m_i$,
which follows from \eqref{04}, \eqref{06} and \eqref{35},
$\hat U^n_i$ are $C$-tight.
In view of \eqref{16} and the fact $\hat X^n_i=\hat U^n_i+\hat L^n_i$, it follows that
$\hat L^n_i$ and $\hat X^n_i$ are also $C$-tight.
\qed

\noi{\bf Proof of Proposition \ref{prop1}.}

i. The $C$-tightness of the sequence follows from Lemma \ref{lem4}
parts (ii) and (iii).

ii.
Consider a convergent subsequence of
$(\hat X^n, \hat L^n,\hat E^n,\hat S^n,\hat P^n,\hat U^n)$, with limit
$(X,L,E,S,P,U)$, where $E$ and $S$ are as before.
Recall \eqref{39}
and note that one has $\hat P^n_i-\hat P^{\#,n}_i\to0$ in probability.
Also note by the second part of \eqref{07} and the tightness
of $\hat L^n_i$, that $T^n_i\to\iota$ in probability.
This and the $C$-tightness of $\hat S^n_i$ implies that
$\hat S^n_i(T^n_i)-\hat S^n_i\to0$ in probability.
Hence, letting $B_i=\sig_i^{-1}(E_i-S_i)$, one has
\[
U_i(t)=X_{0,i}+\sig_iB_i(t)+P_i(t)+m_it.
\]
The map $\Gam$ is continuous in the topology of uniform convergence on compact. Hence, by \eqref{11+}, $(X_i,L_i)=\Gam(U_i)$ for $i\in[N]$, and then by Lemma \ref{lem:sk},
\[
X_i(t)=X_{0,i}+\sig_i B_i(t)+m_it+P_i(t)+L_i(t),
\qquad \int_0^\iy X_i(t)dL_i(t)=0.
\]

Now, $P_i$ is a.s.\ Lipschitz and therefore a.s.\ a.e.\ differentiable. The existence of a progressively measurable process
that is a.e.\ the derivative of $P_i$ follows by an argument that was given in the proof of \cite[Theorem 3.4]{atabudram}.
We denote this derivative by $\beta_i$. 

The goal now is to prove that \eqref{SDIR} is satisfied. By invoking Skorokhod's representation theorem, we may assume without loss of generality that convergence holds a.s. In what follows, fix $\om$ in the full $\PP$-measure set where convergence holds.

Fix a time interval $[0,t]$.
The proof will be complete once it is shown that
\[
\leb(G)=t \qquad \text{where}\qquad
G=\{s\in[0,t]:\beta(s)\in\conv\{b_\pi:\pi\in\Pi(X(s))\}\}.
\]
For $\eps>0$ let
\[
\Pi^\eps(x)=\{\pi\in\Pi:x_i<x_j-4\eps \text{ implies } \pi(i)<\pi(j)\},
\qquad x\in\R^N,
\]
\[
G^\eps=\{s\in[0,t]:\beta(s)\in\conv\{b_\pi:\pi\in\Pi^\eps(X(s))\}\}.
\]
Then $\eps\mapsto\Pi^{\eps}(x)$ is setwise increasing, and
\[
\bigcap_{\eps>0}\Pi^\eps(x)=\Pi(x).
\]
As a consequence, $\eps\mapsto G^\eps$ is setwise increasing and
\[
\bigcap_{\eps>0}G^\eps=G.
\]
By continuity of measure, it suffices to show that for every $\eps>0$,
$\leb(G^\eps)=t$. To this end, we first show the following.
\begin{equation}\label{31}
\begin{split}
&\text{\it For every $\eps>0$ there are $\del_0$ and $n_0$ such that
for $s\in(0,t]$, $n>n_0$ and $\del<\del_0$,}
\\
&\hspace{9em}
\calR^n(\theta)\in\Pi^\eps(X(s)), \qquad \theta\in[s,s+\del].
\end{split}
\end{equation}
To show \eqref{31}, let $n_0$ be so large that for all $n>n_0$,
\[
\max_i\|\hat X^n_i-X_i\|^*_t<\eps.
\]
Let $\del_0>0$ be so small that
\[
\max_i w_t(X_i,\del_0)<\eps.
\]
To show that $\calR^n(\theta)$ is in $\Pi^\eps(X(s))$ is to show that whenever $X_i(s)<X_j(s)-4\eps$, one has $\calR^n_i(\theta)<\calR^n_j(\theta)$, that is, $\rank(i;\hat X^n(\theta))<\rank(j;\hat X^n(\theta))$. The latter will be guaranteed if $\hat X^n_i(\theta)<\hat X^n_j(\theta)$ for all $\theta\in[s,s+\del_0]$. But
\[
\hat X^n_j(\theta)-\hat X^n_i(\theta)> X_j(\theta)-X_i(\theta)-2\eps
> X_j(s)-X_i(s)-4\eps>4\eps-4\eps=0.
\]
This proves \eqref{31}.

In view of \eqref{31}, and recalling $b=\la_0p$, we have for $s\in(0,t]$ and $\del$ and $n$ as above,
\[
\del^{-1}(\hat P^{\#,n}(s+\del)-\hat P^{\#,n}(s))
=\del^{-1}\la_0\int_s^{s+\del}p_{\calR^n(\theta)}d\theta
\in\conv\{b_\pi:\pi\in\Pi^\eps(X(s))\}.
\]
Since for every $x\in\R^N$, $\conv\{b_\pi:\pi\in\Pi^\eps(x)\}$
is a closed subset of $\R^N$, we also have
\[
\del^{-1}(P(s+\del)-P(s))\in\conv\{b_\pi:\pi\in\Pi^\eps(X(s))\}.
\]
For a.e.\ $s$, the limit of the lefthand is $\beta(s)$.
This shows $\leb(G^\eps)=t$ and completes the proof of part (ii).

iii. The tightness shown in part (i), the fact that limits are supported on solutions to \eqref{SDIR}, as shown in part (ii),
and the uniqueness stated in Proposition \ref{prop0}, imply that the entire sequence converges
in distribution to the unique weak solution of \eqref{SDIR}.
\qed

\subsection{The limit under condition \texorpdfstring{\eqref{ICalpha}}{(ICa)}}
\label{sec42}

The goal here is to prove Proposition \ref{prop2},
which is the analogue of Proposition \ref{prop1}.
In this subsection, \eqref{ICalpha} is assumed throughout.

\begin{proposition}\label{prop2}
i.
The sequence $(\check X^n,\hat E^n,\hat S^n,\hat P^{\#,n})$ is $C$-tight, and $\hat L^n\to0$ in probability.
\\
ii. If $(X,E,S,P)$ is a subsequential weak limit, then $(X,B)$ forms a weak solution to \eqref{SDI} with the data indicated in Theorem \ref{th2}, and where $B_i=\sig_i^{-1}(E_i-S_i)$ (with $\sig_i$ as in \eqref{011}) and the progressively measurable rank-dependent drift $\beta$ is the a.e.\ derivative of $P$.
\\
iii. Consequently, denoting $\hat B^n_i=\sig_i^{-1}(\hat E^n_i-\hat S^n_i)$, one has $(\hat X^n,\hat B^n)\To(X,B)$, where the latter is a (weak) solution of \eqref{SDI}.
\end{proposition}

\proof
The arguments are very similar to, only somewhat simpler than those given in Subsection \ref{sec41}. Hence, we only indicate the differences.

Lemma \ref{lem1} holds, and its proof is valid as is, as it has nothing to do with conditions \eqref{IC0} or \eqref{ICalpha}. The same is true with respect to Lemma \ref{lem4} parts i and ii. 

Also, in part iii of Lemma \ref{lem4}, the statements regarding $\hat P^n_i$ and $\hat M^n_i$ and their proof are valid without any change. As for the remaining content of Lemma \ref{lem4} part iii, we prove instead
\begin{equation}\label{102}
\text{
$\hat L^n_i\to0$ in probability, and $\check X^n_i$ are $C$-tight.}
\end{equation}
To this end, subtract $n^{-1/2}\al_n$ from both sides of \eqref{10}
to obtain
\begin{align}\label{10a}
\check X^n_i&=\check U^n_i+\hat L^n_i \qquad \text{where}\qquad
\check U^n_i(t)=\check X^n_i(0-)+\hat E^n_i(t)+\hat A^n_i(t)
-\hat S^n_i(T^n_i(t))+\hat m^n_it.
\end{align}
Now, $\check U^n_i$ are $C$-tight because $\check X^n(0-)$
converge under condition \eqref{ICalpha}, and the remaining terms
in the definition of $\check U^n_i$ are $C$-tight as already shown.

To prove the claim regarding $\hat L^n_i$, note that if $\hat L^n_i(T)>0$ then there exists a time $t\in[0,T]$ such that $X^n_i(t)=0$ and $L^n_i(t)=0$. Therefore, $0=\hat X^n_i(t)=\check X^n_i(t)+n^{-1/2}\al_n$, hence by \eqref{10a}, $\check U^n_i(t)=-n^{-1/2}\al_n$. Thus
\[
\PP(\hat L^n_i(T)>0)\le \PP(\|\check U^n_i\|^*_T\ge n^{-1/2}\al_n)\to0,
\]
where the last statement follows from the tightness of $\|\check U^n_i\|^*_T$ and the fact that $n^{-1/2}\al_n\to\iy$. This proves that $\hat L^n_i\to0$ in probability. In view of \eqref{10a}, this also shows that $\check X^n$ are $C$-tight, and \eqref{102} is proved.

Next, if $(X,U)$ is a limit point of $(\check X^n,\check U^n)$, then we have shown that $X=U$ (compare with $(X_i,L_i)=\Gam(U_i)$ in the case of Subsection \ref{sec41}). 

Based on these statements, the completion of the proof of
Proposition \ref{prop2} follows closely that of Proposition \ref{prop1},
where, in particular, the satisfiability of \eqref{SDI}
is completely analogous to that of \eqref{SDIR}.
\qed

\section{Proof of main results}\label{sec5}

Going back to the steps listed in subsection \ref{sec25}, note that Proposition \ref{prop0} and Remark \ref{rem0} establish steps 1 and 2, and Propositions \ref{prop1} and \ref{prop2} give steps 3 and 4. Steps 5 and 6 are carried out below, which completes
the proof of the main results.

\noi{\bf Proof of Theorem \ref{th2}.} 
The convergence to the solution of \eqref{SDIR} and \eqref{SDI} has already been shown in Proposition \ref{prop1}(iii) and Proposition \ref{prop2}(iii), respectively. To prove the theorem, it remains to show that if $(X,L,B)$ is a weak solution of \eqref{SDIR} then it is also a weak solution of \eqref{SDER}, and similarly, if $(X,B)$ is a weak solution of \eqref{SDI} then it is also a weak solution of \eqref{SDE}. To this end, it suffices to show that, a.s., for a.e.\ $t$, for every $i\ne j$, $X_i(t)\ne X_j(t)$. 

Fix $i \neq j$. For the solution $X$ to \eqref{SDIR}, the difference $X_i(t) - X_j(t)$ is given by  
\[
X_i(t) - X_j (t) = X_{0,i} - X_{0,j} + \sigma_i B_i (t) - \sigma_j B_j (t) + (m_i - m_j) t + \int^t_0 (\beta_i(s) - \beta_j (s) ) d s  + L_i(t) - L_j (t) 
\]
for $t \ge 0 $. By Tanaka's formula, we obtain
\[
\lvert X_i(t) - X_j (t) \rvert = \lvert X_{0,i}- X_{0, j} \rvert + \int^t_0 \text{sgn} ( X_i(s) - X_j (s)) d ( X_i (s) - X_j (s) ) + L^{i,j}(t) 
\]
for $t \ge 0 $, where $\text{sgn} (x) := 1_{\{x > 0 \}} - 1_{\{ x \le 0 \}}$, $x \in \mathbb R$ and $L^{i,j}(\cdot)$ is the local time accumulated at the origin for the semimartingale $X_i(\cdot) - X_j (\cdot)$. 
Take a nonnegative function $\phi \in C^2_b (\mathbb R_+) $ with the nonincreasing, nonnegative second derivative $\phi^{\prime\prime}$ that satisfies $\phi^{\prime\prime} (u) = 1 $ for $u \in [0, 1/2]$, $\phi^{\prime\prime}(u) = 0 $ for $u \ge 1$ and $\phi(u) = \phi^\prime (u) = 0 $ for $u \ge 1$. Then applying It\^o's formula to $\phi ( \varepsilon^{-1} \lvert X_i (t) - X_j (t) \rvert ) $, we obtain 
\[
\hspace{-5cm} \phi \Big( \frac{1}{\varepsilon} \lvert X_i (t) - X_j (t) \rvert \Big) - \phi \Big( \frac{1}{\varepsilon} \lvert X_{0,i} - X_{0,j} \rvert \Big) 
\]
\[
= \frac{1}{\varepsilon} \int^t_0 \phi^\prime \Big(\frac{1}{\varepsilon} \lvert X_i (u) - X_j (u) \rvert \Big) d \lvert X_i (u) - X_j (u) \rvert + \frac{1}{2\varepsilon^2} \int^t_0 \phi^{\prime\prime} \Big ( \frac{1}{\varepsilon} \lvert X_i (u) - X_j (u) \rvert \Big) (\sigma_i^2 + \sigma_j^2) d u 
\]
for $\varepsilon > 0 $ and $t \ge 0 $. This implies that
\[
\int^t_0 \phi^{\prime\prime} \Big ( \frac{1}{\varepsilon} \lvert X_i (u) - X_j (u) \rvert \Big) d u = \frac{2\varepsilon^2}{(\sigma_i^2 + \sigma_j^2) } \Big( \phi \Big( \frac{1}{\varepsilon} \lvert X_i (t) - X_j (t) \rvert \Big) - \phi \Big( \frac{1}{\varepsilon} \lvert X_{0,i} - X_{0,j} \rvert \Big) \Big)  
\]
\[
\hspace{4cm} - \frac{2\varepsilon}{\sigma_i^2 + \sigma_j^2}\int^t_0 \phi^\prime \Big(\frac{1}{\varepsilon} \lvert X_i (u) - X_j (u) \rvert \Big) d \lvert X_i (u) - X_j (u) \rvert 
\]
for $\varepsilon > 0 $ and $t \ge 0 $. Taking the limits as $\varepsilon \downarrow  0$, we obtain 
\[
\liminf_{\varepsilon \downarrow 0 } \int^t_0 \phi^{\prime\prime} \Big ( \frac{1}{\varepsilon} \lvert X_i (u) - X_j (u) \rvert \Big) d u = 0 
\]
almost surely. Since $\varphi^{\prime\prime}(0) = 1 $, this implies, by Fatou's lemma, that
\begin{equation}
\begin{split}
\int^t_0 1_{\{ X_i (u) = X_j (u) \}} d u &= \int^t_0 \liminf_{\varepsilon \downarrow 0 } \phi^{\prime\prime} \Big ( \frac{1}{\varepsilon} \lvert X_i (u) - X_j (u) \rvert \Big) d u \\
& \le \liminf_{\varepsilon \downarrow 0 } \int^t_0 \phi^{\prime\prime} \Big ( \frac{1}{\varepsilon} \lvert X_i (u) - X_j (u) \rvert \Big) d u = 0
\end{split}
\end{equation}
Therefore,  the set of times when two or more components $X_i$ meet
is shown to have Lesbegue measure zero. If $(X, L, B)$ is a weak solution of  \eqref{SDIR}, then it is also a weak solution of \eqref{SDER}. 

The solution of \eqref{SDI} can be handled in a similar manner. If $(X, B)$ is a weak solution of \eqref{SDI}, then it is also a weak solution of \eqref{SDE}. 
\qed

\skp

\noi{\bf Proof of Theorem \ref{th1}.} 
Pathwise uniqueness for \eqref{SDI}, \eqref{SDIR}, \eqref{SDE} and \eqref{SDER} has been shown in Proposition \ref{prop0} and Remark \ref{rem0}. Weak existence for the four equations follows from Theorem \ref{th2}. It remains to prove strong existence.
To this end we employ the Yamada-Watanabe Theorem (see Remark \ref{rem: Yamada-Watanabe-Kurtz}).
\qed

\skp

\noi
{\bf Acknowledgement.}
The first author is partially supported by ISF grant 1035/20. The second author is partially supported by NSF grant DMS-2008427. The authors would like to thank the Isaac Newton Institute for Mathematical Sciences, Cambridge, for support and hospitality during the programme Stochastic systems for anomalous diffusions (SSD), where a part of the work on this paper was undertaken.


\bibliographystyle{is-abbrv}

\bibliography{main}

\vspace{.5em}

\end{document}